\documentclass[12pt]{article}
\usepackage{amsmath}
\usepackage{graphicx}
\usepackage{subfigure}

\usepackage[affil-it]{authblk}
\usepackage[english]{babel}
\usepackage{blindtext}

\begin{document}

\title{Mutual information for fitting deep nonlinear models}
\author[1]{Jacob S. Hunter}
\author[1]{ Nathan O. Hodas \thanks{Corresponding author}}
\affil[1]{ Pacific Northwest National Laboratory}
\affil[ ]{\textit {\{jacob.hunter, nathan.hodas\}@pnnl.gov}}
\date{\today}
\maketitle

\begin{abstract}
Deep nonlinear models pose a challenge for fitting parameters due to lack of knowledge of the hidden layer and the potentially non-affine relation of the initial and observed layers.  In the present work we investigate the use of information theoretic measures such as mutual information and Kullback-Leibler (KL) divergence as objective functions for fitting such models without knowledge of the hidden layer.  We investigate one model as a proof of concept and one application of cogntive performance.  We further investigate the use of optimizers with these methods.  Mutual information is largely successful as an objective, depending on the parameters.  KL divergence is found to be similarly succesful, given some knowledge of the statistics of the hidden layer.
\end{abstract}

\section{Introduction}

There is extensive literature on the effects of sustained activity on human cognitive performance, which suggests that performance is determined by some finite resource that causes a decline in performance when depleted\cite{Boksem:2005, Boksem:2008, Shah:2012}.  However, the nature of this resource remains poorly understood.  In other work (to be submitted) we propose a model in which standardized test performance depends on cognitive resources that are depleted as questions are answered. In this work, we present a method for fitting ODE's to data when the mapping between the observables and underlying hidden variables is unknown.  For the cognitive performance example, a researcher may propose a microscopic mechanistic model for resource variation in the brain but not know how the variations in the physiological resources translate into variation in the observed performance.

We can characterize the relation of cognitive performance to physiological resource usage as a deep nonlinear system.  We consider a nonlinear system to be deep if it can be represented as a composition of multiple nonlinear transformations.  There is an initial hidden layer, followed by a transformation layer, and an output layer that corresponds to what is actually observed. In short, complex nonlinear dynamics on the microscopic scale ultimately translate to emergent behavior on the macroscale via multiple scales, spatially and temporally.  Many types of objective functions may be used to fit equations to hidden variables, such as Kalman filters, Bayesian models, and simple root-mean-squared error (RMSE) regression via supervised learning\cite{nelles2013nonlinear, peterka1981bayesian, voss2004nonlinear}.  In general, these models are parametric, and their success requires some fundamental assumptions about the dynamics producing data.  For example, to minimize RMSE, one needs to be able to regress to the final observed variable, so the hidden layer needs to also be known or learned. 

One approach to understanding a deep nonlinear system would be to propose a model for the hidden variables, construct a function approximator (such as a neural network) as the transformation layer, and then use a fitting procedure such as backpropagation\cite{nelles2013nonlinear} to jointly regress the generic transformation layer and fit the underlying hidden model. The viability of this approach depends on the fidelity of  the function approximator, which may require an extremely large amount of data to properly determine, akin to traditional deep learning.  Here, we propose a method that is agnostic to the transformation layer by seeking to maximize explanatory power of the underlying hidden model to the ultimate observables.

We investigate information theoretic measures, such as mutual information and Kullback-Leibler divergence, as objective functions for fitting the model to generated data.  We construct a number of example nonlinear systems to validate the ability for these information theoretic measures to `peer through' the intervening nonlinear unknown layers and directly fit observations to the hidden variables we are attempting to model.   We confirm that mutual information may be a valuable tool for fitting  models to deep nonlinear systems.

\subsection{Deep Nonlinear Models}

In the present work, we consider systems that we may call ``deep nonlinear models," in that we may have a clear hypothesis for what the dynamics of the initial layer may be,  but there remains an unknown layer that transforms the output of the dynamics into what is observed.  We may consider deep nonlinear models to take the form
\begin{align}
\dot{x} &= f(x,t) &\text{Initial Layer} \\
\dot{y} &= g(x,y,t) &\text{Hidden Layer} \\ 
z &= h(y,t). &\text{Observable/Measurement Layer}
\end{align}
where $f,g,h$ are arbitrary transformations. 
 
This hidden layer may be due to some unknown dynamics, stochasticity, etc.  For example, consider the task of using the number of people at the beach to infer the weather. We may know how to predict the weather (initial layer), and we may separately have an empirical count of beach-goers (observation), but the transformation from the weather to how people decide when to visit the beach remains hidden.  Because we only care about making estimates of the weather, we would prefer to not have to model people's decision making processes explicitly.

Frequently, we may only have the observations -- and not even know the initial layer.  For example, we may wonder if we can count  people at the beach to estimate the weather at different locations.  More generally, we may attempt some sort of model identification given some limited sensor activity.  To be successful, we need to `see through' the hidden layers to directly fit the parameters of the initial layer.  Many techniques exist to do this.  Here, we present an information theoretic technique based on mutual information that allows us to non-parametrically conduct fitting under weak assumptions about the dynamics of the hidden layer.

\subsection{Mutual Information}
Mutual information ($MI$) quantifies how many bits of entropy we may reduce our uncertainty in  $X$  by knowing another variable $Y$ and vice-versa. For example, how much information about an individual's performance ($X$) do we obtain by knowing the amount of cognitive resources available to him or her ($Y$)? $MI$ is defined as $MI(X,Y) = H(X) + H(Y) - H(XY)$, where $H(X)$ is the entropy of the random variable $X$, and $H(XY)$ is the entropy of the joint distribution of $X$ and $Y$.  In the present work, we leverage the property that $MI(X,Y) = MI(f(X),g(Y))$ for invertible functions $f$ and $g$~\cite{Nair:2006ws}, so the underlying initial layer of the nonlinear system may be fit from the observations as long as the two can be related via an isomorphism.  If the relations between $X$ and $Y$ are not isomorphic, we may still utilize mutual information, but we will only be able to explain limited portions of the dynamics. In contrast, optimizing a regression model or Pearson correlation for  parameter estimation requires not only correctly modeling the initial layer but also knowing how those resources quantitatively translate into the observations, because $R^2$ is only maximized when predicted and observed results have an affine relation.

\begin{figure}[htbp]
   \centering
   \includegraphics[width=0.48\columnwidth]{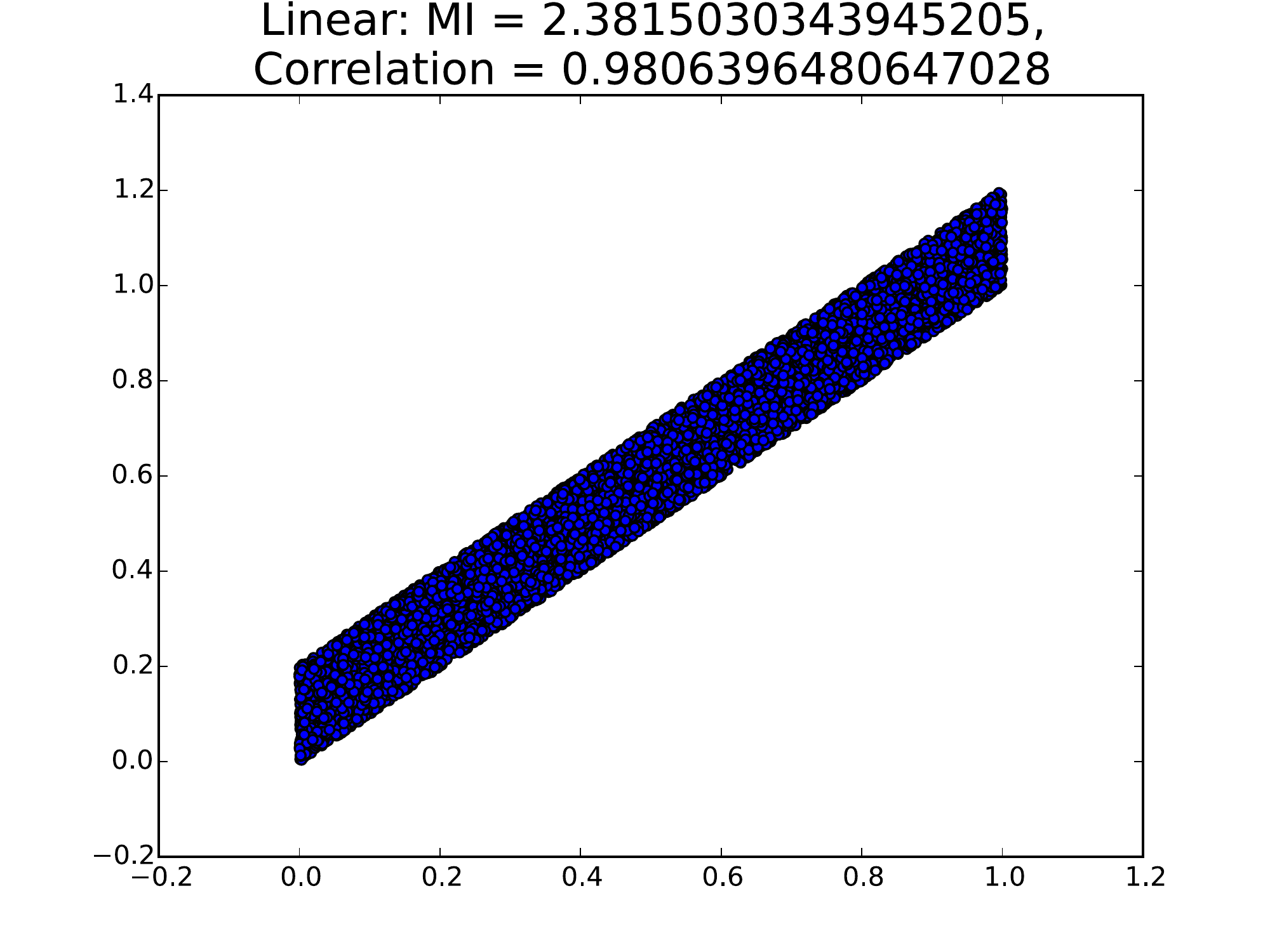} 
   \includegraphics[width=0.48\columnwidth]{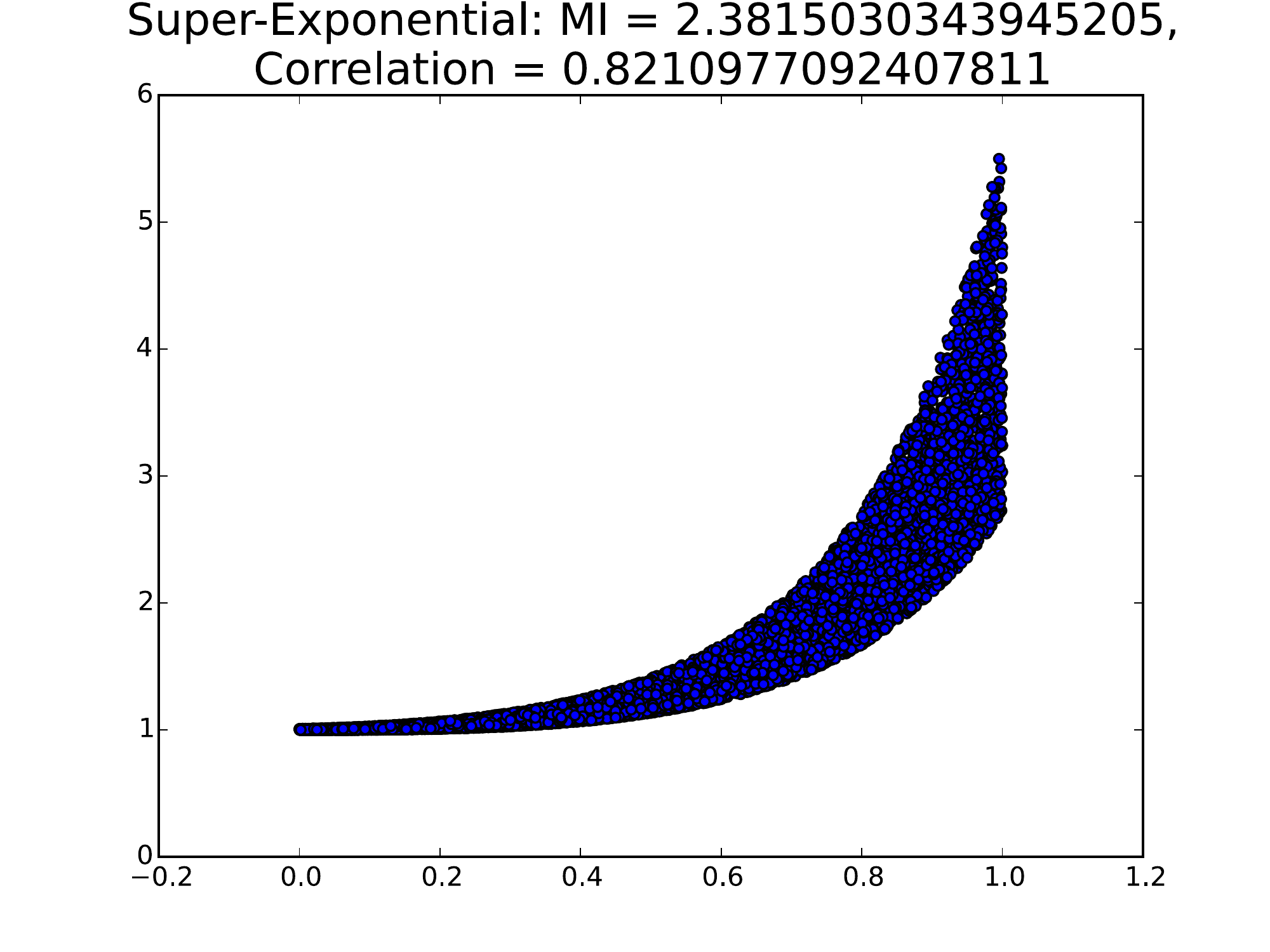} 
   \includegraphics[width=0.48\columnwidth]{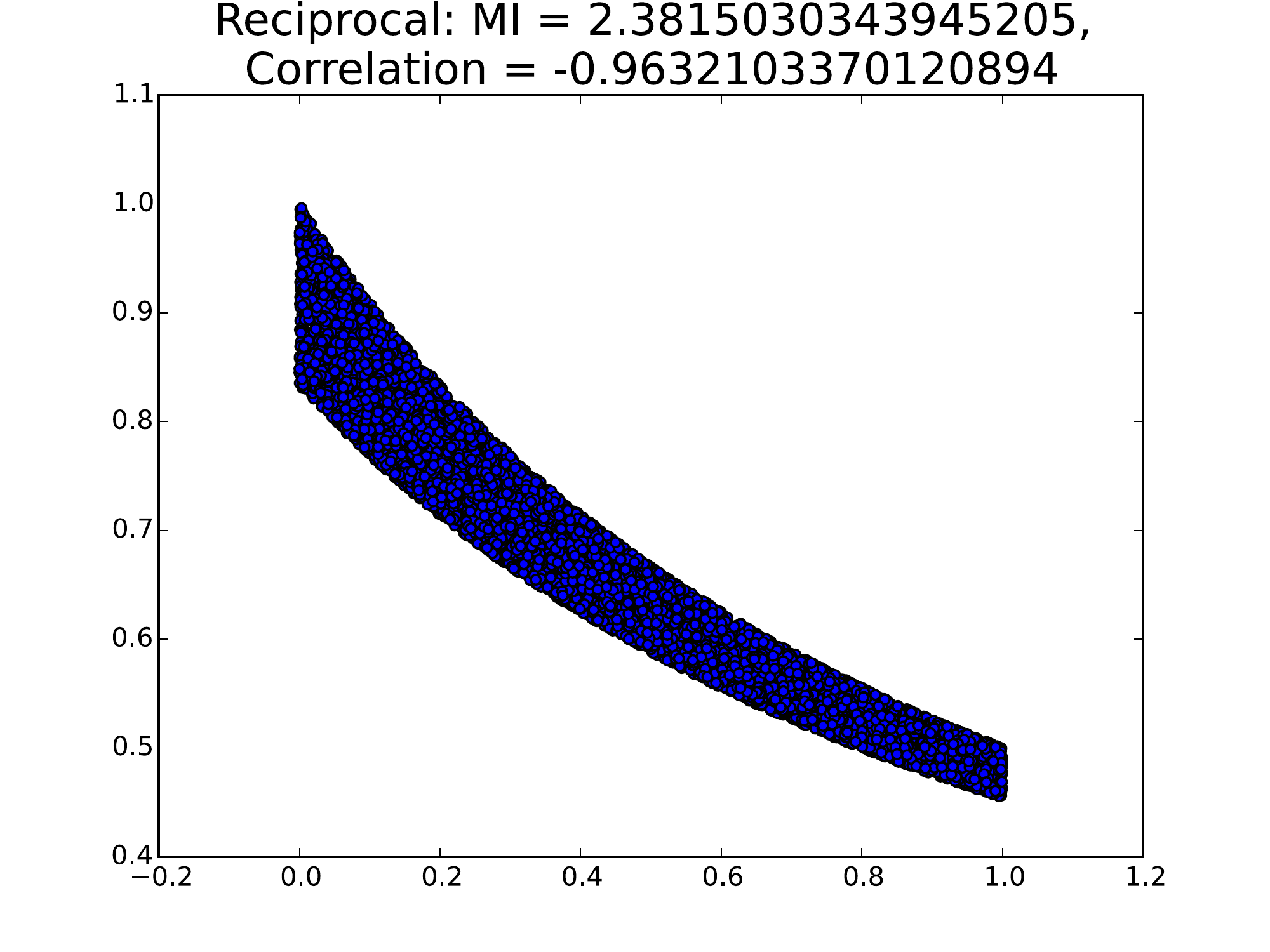} 
   \includegraphics[width=0.48\columnwidth]{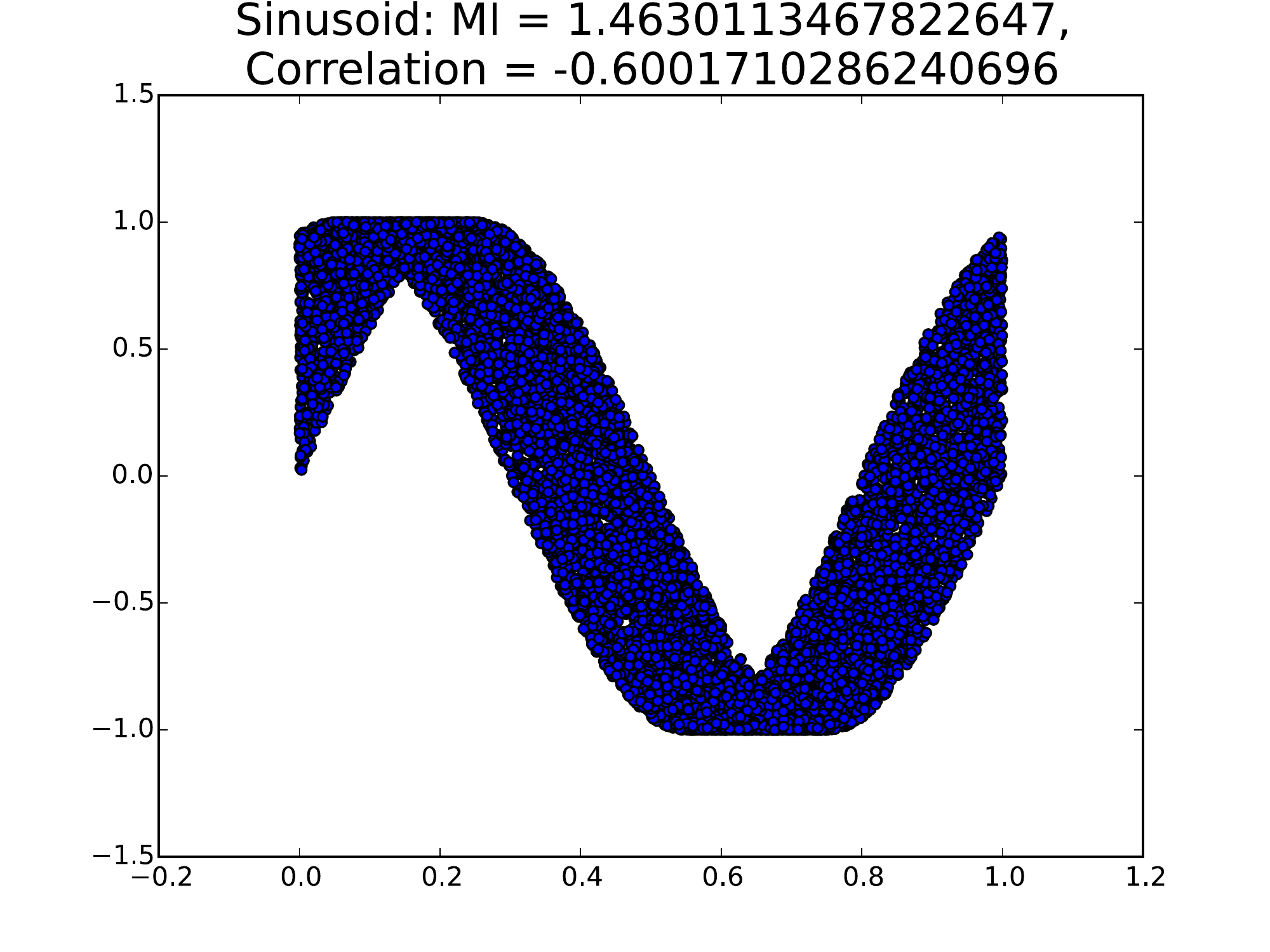} 
   \caption{Mutual information does not vary under invertible transformations, while correlation is only maximize with a perfect linear relationship.}
   \label{fig:mi}
\end{figure}

For example, in Fig.~\ref{fig:mi}, we demonstrate how the mutual information between two variables remains unchanged under isomorphisms. However,  the correlation between the transformed variables may go from positive to zero to negative, depending on the transformation.  Only the cosine transformation, which is not an isomorphism, shows a decrease in mutual information.

This ability to quantitatively relate the explanatory power of two variables, regardless of the isomorphic transformation between them, allows this method to relate the dynamics of the initial layer with the observation in the final layer.  In addition, mutual information is well defined for discrete, continuous, or categorical data -- or some combination thereof.  In this work, we present a method for using mutual information to fit dynamics of the initial layer to discrete observations, without  having to know the intervening hidden layer.

\subsection{KL Divergence}
In addition to mutual information, we tested the Kullback-Leibler (KL) divergence as an objective function.  KL divergence is defined as the amount of information lost when one distribution is approximated by a second; it can thus be thought of as a distance between the two distributions, although it is not symmetric.  The  KL divergence of continuous random variables from $Q$ to $P$ is given by
\begin{equation}
D_{KL}(P||Q) = \int_{-\infty}^{\infty}p(x)\ln\frac{p(x)}{q(x)} dx
\end{equation}
When testing KL divergence as an objective function, we considered the distribution of observations against our initial model under a variety of comparisons.  First, we considered the divergence from the fitted model to the resource values used to generate the outcome data, given the outcome.  In this case, the objective is a loss function, as the goal is for the model to closely approximate reality.  Finally, we considered the divergence from the distribution of resource values given one outcome to the distribution given the other outcome; when this divergence is maximized, it forces the distributions to be as different as possible, which should increase the explanatory power of the model.  

\section{Approach}

\subsection{Example with Mutual Information}
As an example and proof of concept of our fitting with mutual information, we construct a family of deep nonlinear models:
\begin{align}
x &= e^{-\lambda t}  \\
y &= x +y_ \text{err}\\
z &= f(y),
\end{align}
with $y_ \text{err}$ being Gaussian noise and $z = f(y)$ having one of three functional forms:
\begin{align}
f_1(y)& = a y\\
f_2(y) &= e^{ay}\\
f_3(y) &= \sin (a y).
\end{align}
In each case, we wish to fit $\lambda$ to the generated data given the time series of $x$ and $z$.  To generate the data, we chose a value of $\lambda$ and a  range for $t$.  With a fixed $\lambda$ and a given $t$, we can then easily determine $f_1$, $f_2$, and $f_3$.  

To evaluate if a maximum in  mutual information occurs with the true $\lambda$, we simply plotted the $MI$ between data generated from $\lambda$ and different values of $x$ generated from a series of values of guesses $\hat{\lambda}$:  we  varied $\hat{\lambda}$ in increments of 0.01 in the range (-1, 4) and calculated the mutual information $MI(z, e^{-\hat{\lambda} x})$.  This entire process was repeated for several values of $\lambda$ and $a$.  Because of the highly deterministic relationship between $x$, $y$ and $z$, we estimated mutual information using the Non-parametric Entropy Estimation Toolbox (NPEET) with Local Non-uniformity Correction (LNC)\footnote{https://github.com/BiuBiuBiLL/NPEET\_LNC}~\cite{Gao:2014, Kraskov:2004}. This produced far more robust, less biased estimates of mutual information, by correcting for systematic errors that can emerge for deterministic data for the Kraskov estimator.%

\begin{figure}[htbp]
   \centering
   \begin{subfigure}[{$a=1,\lambda=0.5$}\label{fig:a1}] {
   \includegraphics[width=0.48\columnwidth]{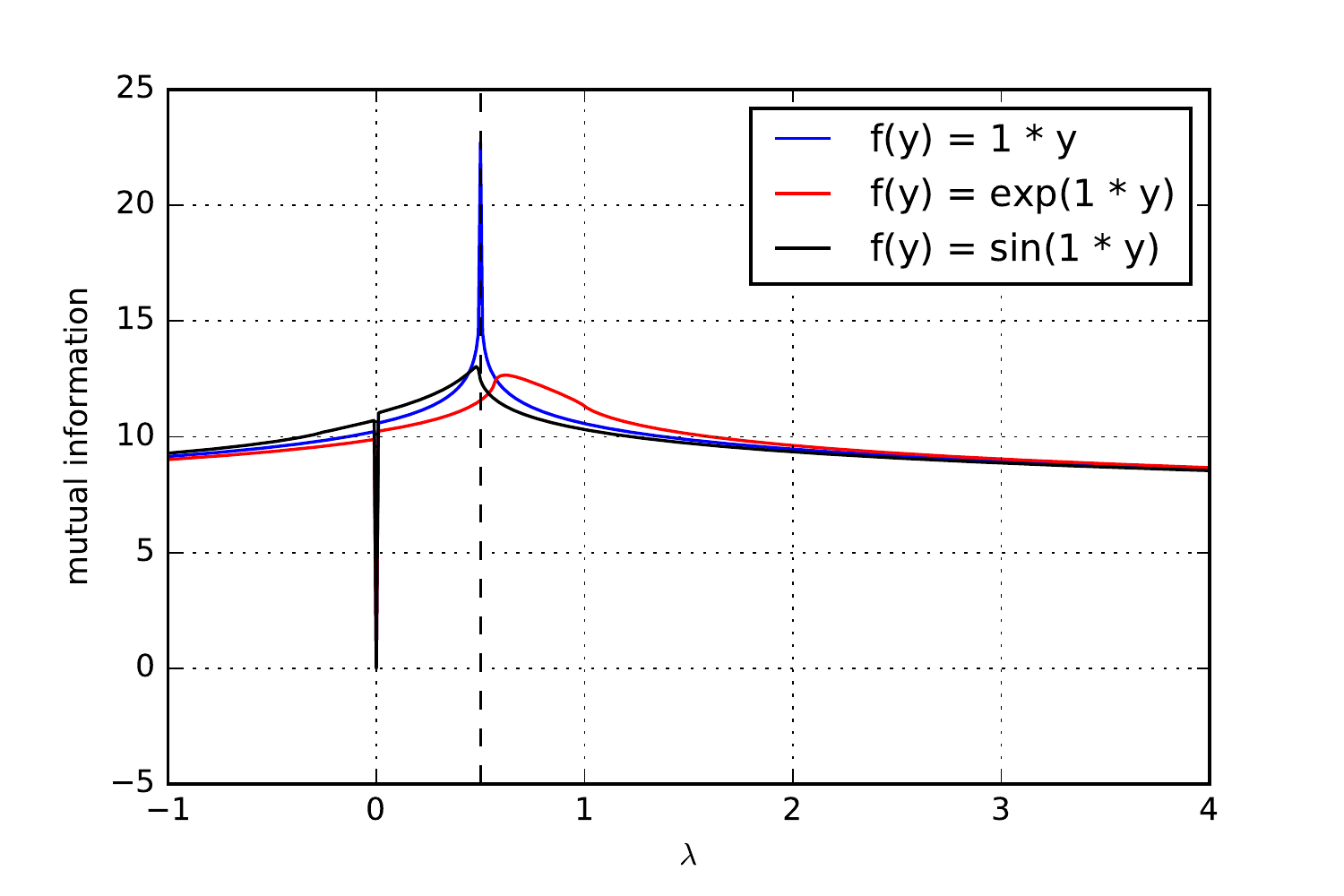} }
   \end{subfigure}
   \begin{subfigure}[{$a=2,\lambda=1$}\label{fig:a2}]{
      \includegraphics[width=0.48\columnwidth]{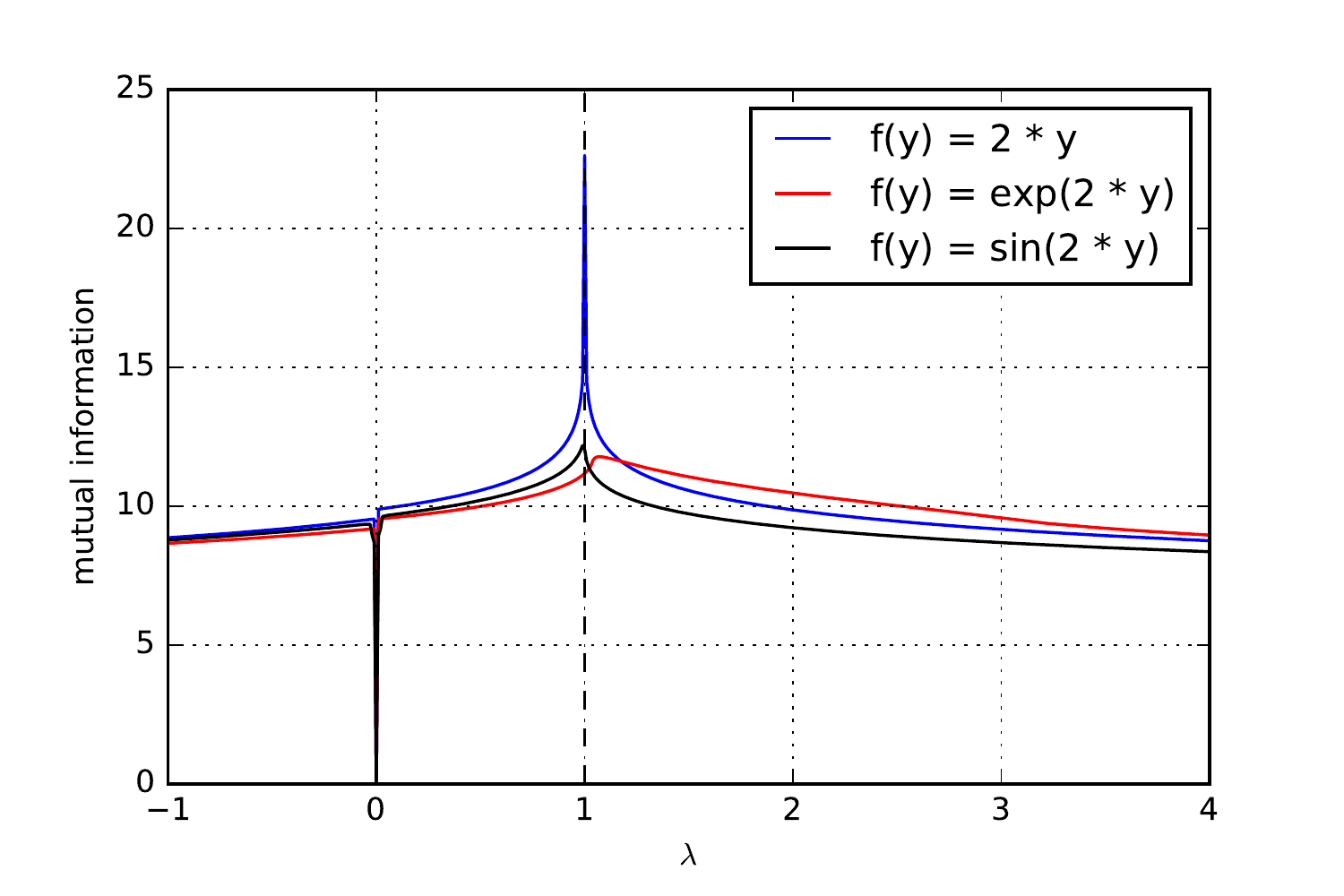} }
   \end{subfigure}
   \begin{subfigure}[{$a=3,\lambda=2$}\label{fig:a3}]{
   \includegraphics[width=0.48\columnwidth]{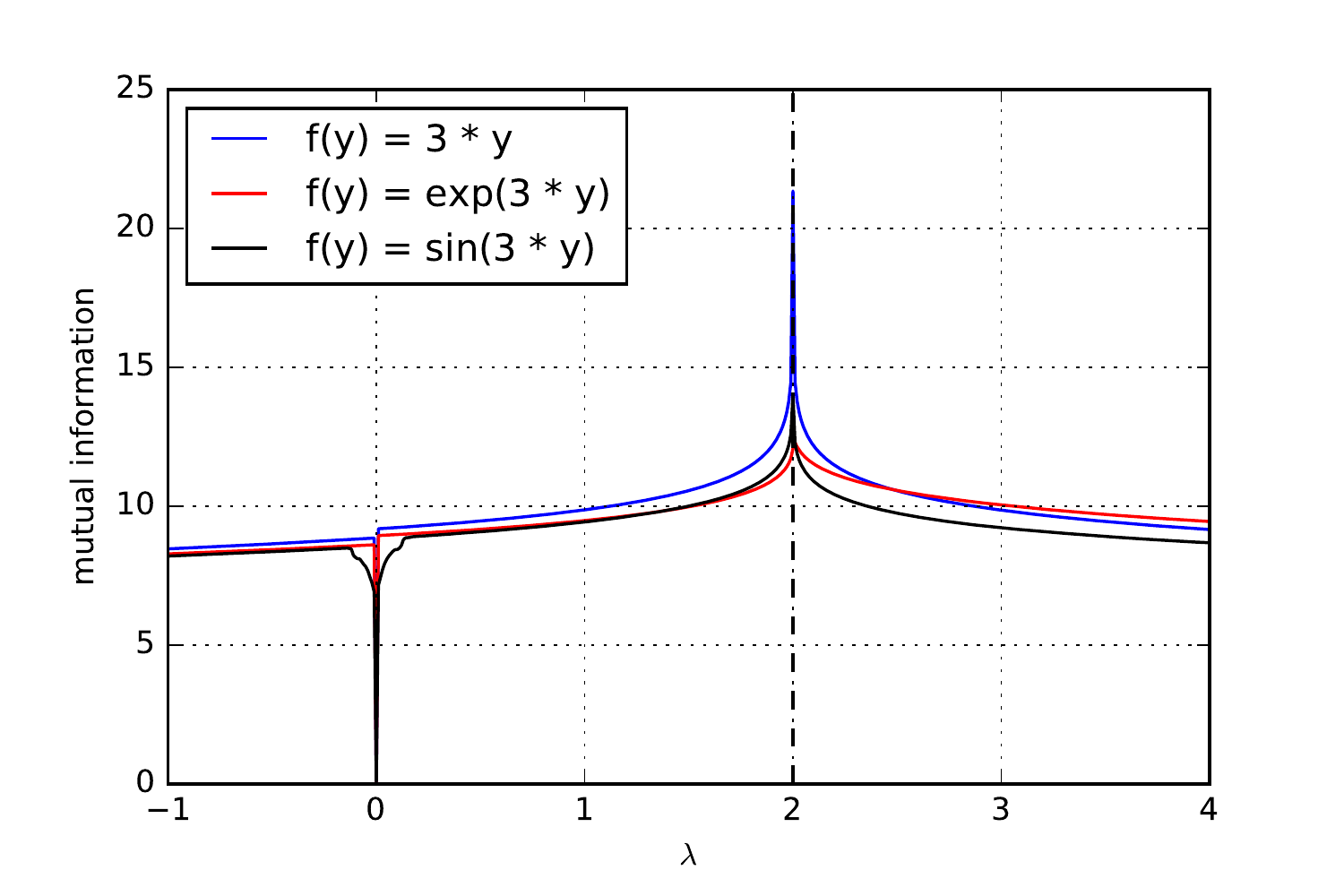} }
   \end{subfigure}
   \caption{Mutual information of outer layer $z=f(y)$ with fitted hidden layer $y'=e^{-\lambda x}$ for different forms of $f(y)$ and true values of $\lambda$.  True value for $\lambda$ is given by dashed line.}
   \label{fig:examplefit}
\end{figure}

Fitting the example model using mutual information correctly produced a maximum of the objective function at or very close to the true value for $\lambda$ (Fig.~\ref{fig:examplefit}).  In particular, Fig.~\ref{fig:a3} shows  a sharp maximum at the true value of $\lambda = 2$ for all forms of $z$.  While the mutual information curves for the exponential and sinusoidal $z$ in Figs.~\ref{fig:a1}~and~\ref{fig:a2} show smoother maxima that are not perfectly at the true value, they provide a fairly accurate estimate of the true value.  The dip in mutual information to 0 at $\lambda = 0$ is due to $y$ taking the form $e^{0x}=1$; in practice, this would also prevent a local optimizer from transitioning over regimes when a particular parameter  completely decouples the hidden model from the observations.  It would require a more global optimization strategy.

\subsection{Application: Cognitive Performance}

A particular system we consider in the present work is cognitive performance over time.  Previous work established a model for cognitive performance based on lactate metabolism~\cite{Cloutier:2009, Wyss:2011, Xu:2007}.  In this model, performance is determined by a primary resource, $A$, which is low during normal resting conditions.  Work on a task consumes $A$.  However, $A$ is replenished by conversion of a secondary resource, $B$, into $A$; $B$ then recovers during resting conditions.  These dynamics are modeled as:
\begin{subequations}\label{eq:cogdep}
\begin{align}
\dot{A}(t) &= w(A,B,t) - \frac{k_w}{t^\rho}\frac{A(t)}{K_A + A(t)}\\
\dot{B}(t) &= -w(A,B,t) + \frac{k_r}{t^\rho}\left(B_{max} - B(t)\right)(1 - \delta(t))\\
w(A,B,t) &= \frac{k_b}{t^\rho}\frac{\left(1 - A\left(t\right)\right)B(t)}{K_B + B(t)} \delta(t),
\end{align}
\end{subequations}
where $w$ is the rate of conversion of $B$ into $A$ and $\delta(t)$ is given by 1 when on task (working) and 0 when off task (resting).  

Although this model may be difficult to understand outside the context of modeling cognitive performance, it captures many of the properties that make fitting difficult for many methods. Clearly, this model is complex, nonlinear and time-heterogeneous.  In addition, it is non-differentiable at the boundaries between  tasks. Furthermore, although this may be a proposed model of cognitive resources, we have no hypothesis for how cognitive resources translate into correct/incorrect answers on a standardized test.  All we can assume is that fewer resources decreases the likelihood of correct answers. Thus, to fit the model, we need to be ambivalent to the intervening hidden layer.

\subsection{Generated Performance Data}
To test fitting the model, we generated data similar to data from online practice standardized tests, obtained from the Kaggle  ``What do you know" competition\footnote{http://www.kaggle.com/c/WhatDoYouKnow} and provided by grockit.com.  Problem times were generated stochastically from an exponential distribution,
\begin{equation}
p(x) = \lambda e^{-\lambda x},
\end{equation}
with $x$ the time spent on or off task, in minutes.  For time on task, we set $\lambda = \frac{1}{4}$; for time off task, $\lambda = \frac{1}{4}$, except for every tenth task, when $\lambda = \frac{1}{40}$.  Our fictional participant thus spent an average of four minutes on or off each task, with a longer break every ten tasks.  To generate the dataset, 25 time series with 3000 event times were generated in this fashion.  From the generated time series, we used a numerical ODE solver from the scipy.integrate library to get the primary and secondary resources from the model using the previously obtained parameters.  Outcome of a task was determined randomly based on the primary resource, $A$, at the end of each question; the probability of success was given by the logistic sigmoid function,
\begin{equation}
P(\text{correct}|A) = (1+e^{-\alpha(A-A_0)})^{-1},
\end{equation}
with $\alpha$ and $A_0$ chosen as 169 and 0.204, respectively, such that typical high values of $A$ would give a success rate of 70\% and typical low values whould give a success rate of 30\%.

Using this generated data, we sought to refit the resource model to determine whether our fitting methods are appropriate.  To start, we fit each parameter individually to get a sense for the shape of the objective function.  To fit a parameter, we varied around the true value and solved the model ODEs for the resources at each time.  We then used the newly calculated model resources and the outcome data to calculate the objective function at each parameter step.

\subsection{Objective Functions}

To construct an objective function using mutual information, we estimated mixed continuous-discrete mutual information using the implementation based on Kraskov et. al. and Kozachenko and Leonenko in NPEET\footnote{http://www.isi.edu/\~{}gregv/npeet.html} \cite{Kraskov:2004, Khan:2007, Kozachenko:1987}.  We calculated the mutual information of the series of task outcomes (with 1 standing for correct/success and 0 for incorrect/failure) with the primary resource at the end of each task.  

KL divergence was also estimated using the implementation in NPEET \cite{Wang:2009, Khan:2007}.  We used two methods for using KL divergence as an objective function.  In both cases, we considered the frequency distribution of the primary resource at the end of a task given the outcome.  
The first use of KL divergence relied upon prior knowledge of the underlying resources used to generate the data.  We again used the distribution of fitted resources given outcome, but now compared to the distribution of the underlying resources: $$D_{KL}\left(p(A_{\text{fitted}}|\text{correct}) || p(A_\text{underlying}|\text{correct})\right)$$ and  $$D_{KL}\left(p(A_\text{fitted}|\text{incorrect}) || p(A_\text{underlying}|\text{incorrect})\right).$$  

In the second method, we compared the frequency distribution of the primary resource given a successful outcome with that given an unsuccessful outcome: $$D_{KL}(p(A|\text{correct}) || p(A|\text{incorrect}))$$ and  $$D_{KL}(p(A|\text{incorrect}) || p(A|\text{correct})).$$  By maximizing either of these divergences or a combination thereof, the model should be as distinct as possible within a regime of cognitive depletion compared to a regime of plentiful cognitive resources.  

\subsection{Optimizers}

Because mutual information does not require the random variables to be positively correlated for it to be maximized, using it as a metric can result in model resource values that are lower on average for successful outcomes than for unsuccessful outcomes.  While this could be interpreted as the model representing waste products instead of resources, it goes against the assumptions that base the model on lactate metabolism.  Thus, we added a constraint that the average resource given success be greater than that given failure:
\begin{equation}
f_c = \langle A|\text{correct}\rangle - \langle A|\text{incorrect}\rangle \ge 0,
\end{equation}
where $A$ is the primary resource modeled in Eq.~\ref{eq:cogdep}.

We ran tests on optimizing all parameters together over the mutual information objective using constrained optimizers implemented in NLOPT\footnote{Steven G. Johnson, The NLopt nonlinear-optimization package, http://ab-initio.mit.edu/nlopt}.  However, the optimizers here posed problems for this task.  Because the objective function requires solving the model ODEs for the time series and parameters, it takes a long time to evaluate.  Additionally, because the  Kraskov mutual information is non-differentiable, the gradient is unknown and must be estimated; therefore, each iteration of the optimizers required multiple evaluations of a computationally expensive objective.  Thus, we use Simultaneouse Perturbation Stochastic Approximation (SPSA).  SPSA lends itself to  objective functions requiring a simulation or other computationally expensive task by requiring only two calls of the objective function to estimate the gradient\cite{Spall:1998, Wang:Spall:2008}.   

\section{Discussion}

\begin{figure}[htbp]
   \centering
   \includegraphics[width=0.48\columnwidth]{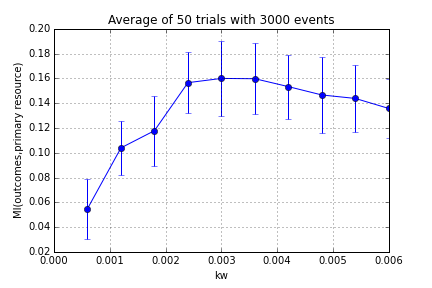} 
   \includegraphics[width=0.48\columnwidth]{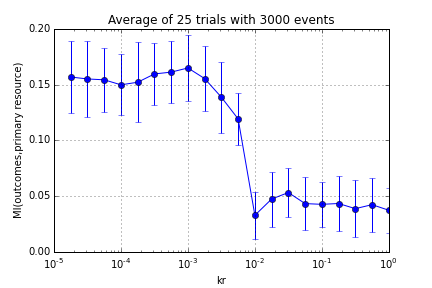}
   \includegraphics[width=0.48\columnwidth]{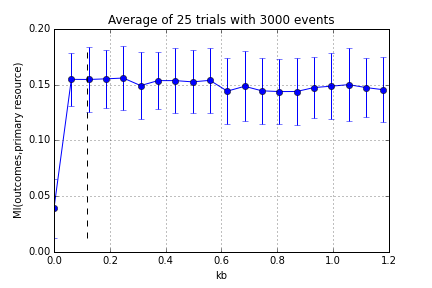} 
   \includegraphics[width=0.48\columnwidth]{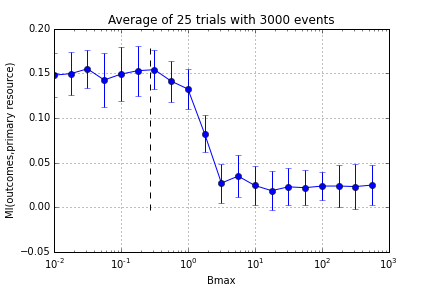} 
   \includegraphics[width=0.48\columnwidth]{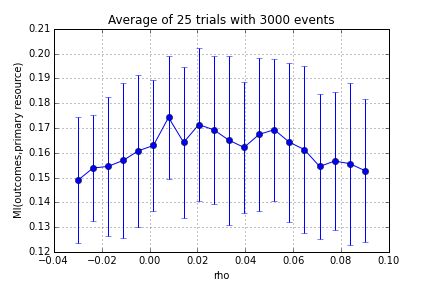} 
   \caption{Mutual information of primary resource with task outcomes.  True parameters are given by dashed lines.}
   \label{fig:mifit}
\end{figure}

When fitting one parameter from the true value and holding the rest constant, mutual information as an objective function produced a maximum at or near the correct parameter values for $k_w$, $k_r$, and $B_{max}$; the landscapes for $\rho$ and $k_b$ featured a plateau-like range of values at the respective maxima (Fig. \ref{fig:mifit}).

\begin{figure}[htbp]
   \centering
   \includegraphics[width=0.48\columnwidth]{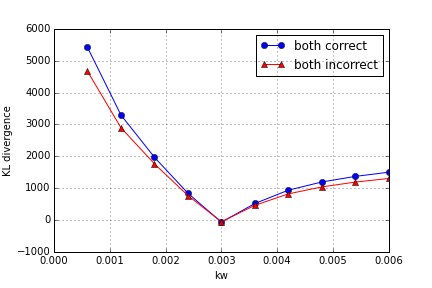} 
   \includegraphics[width=0.48\columnwidth]{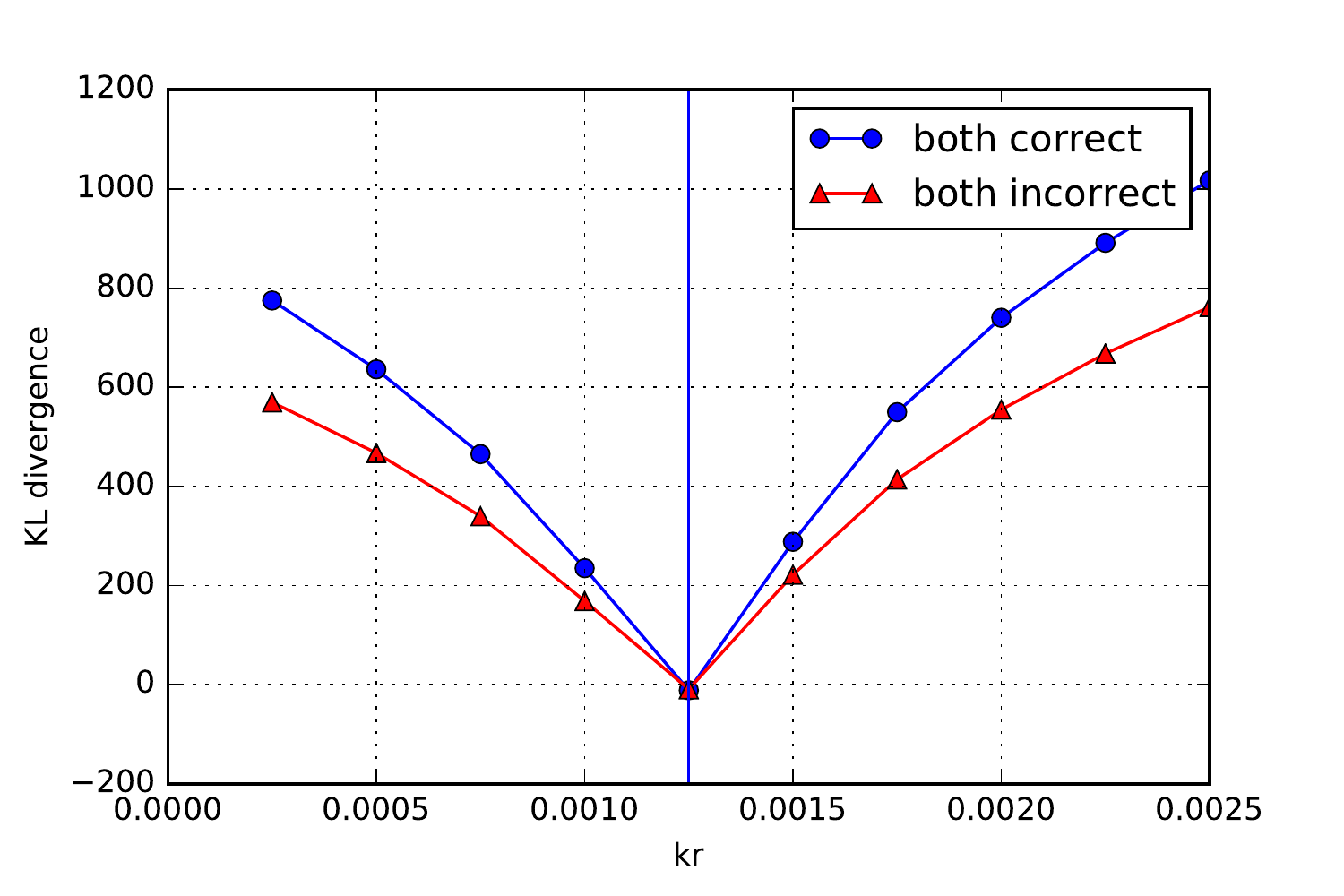} 
   \includegraphics[width=0.48\columnwidth]{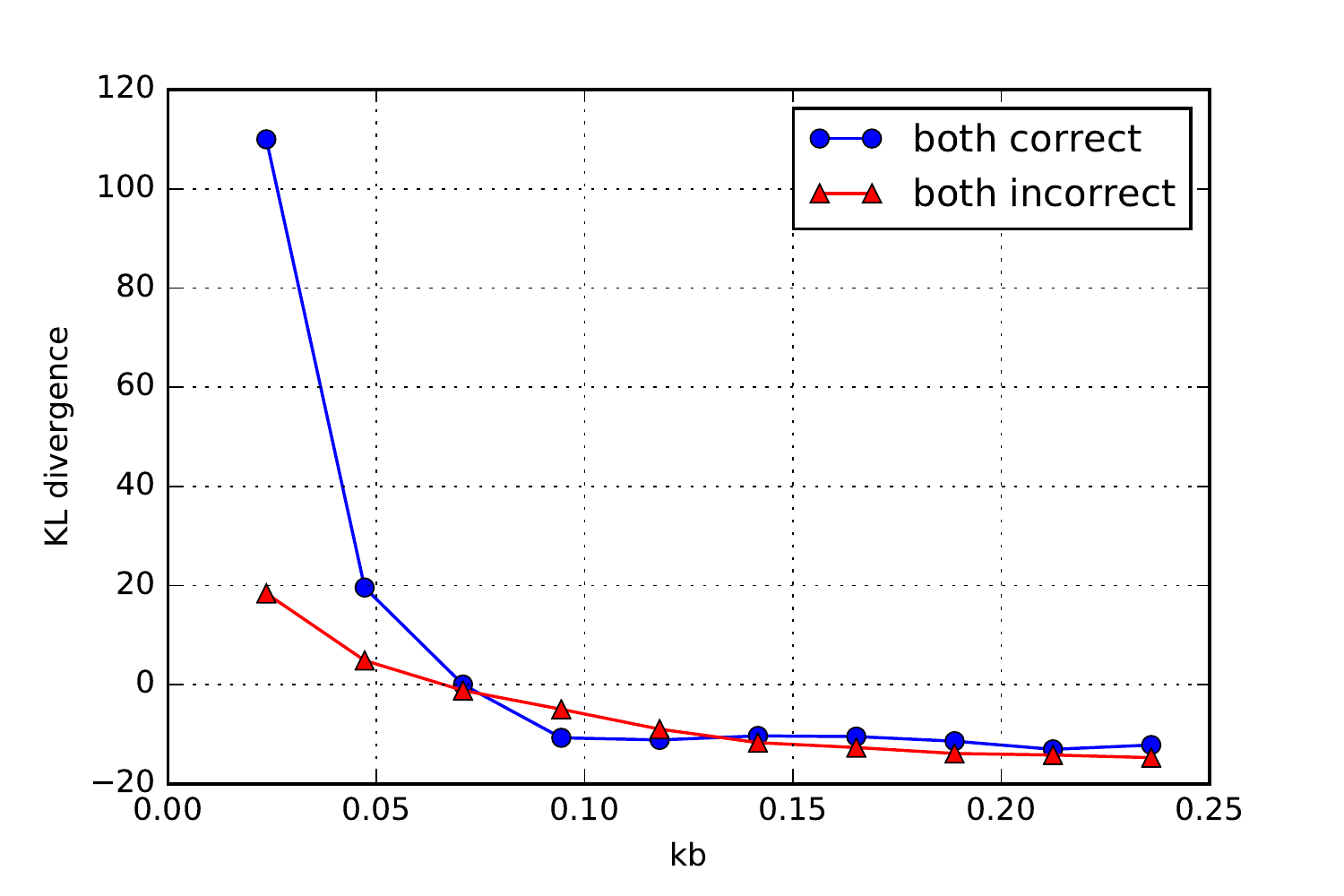} 
   \includegraphics[width=0.48\columnwidth]{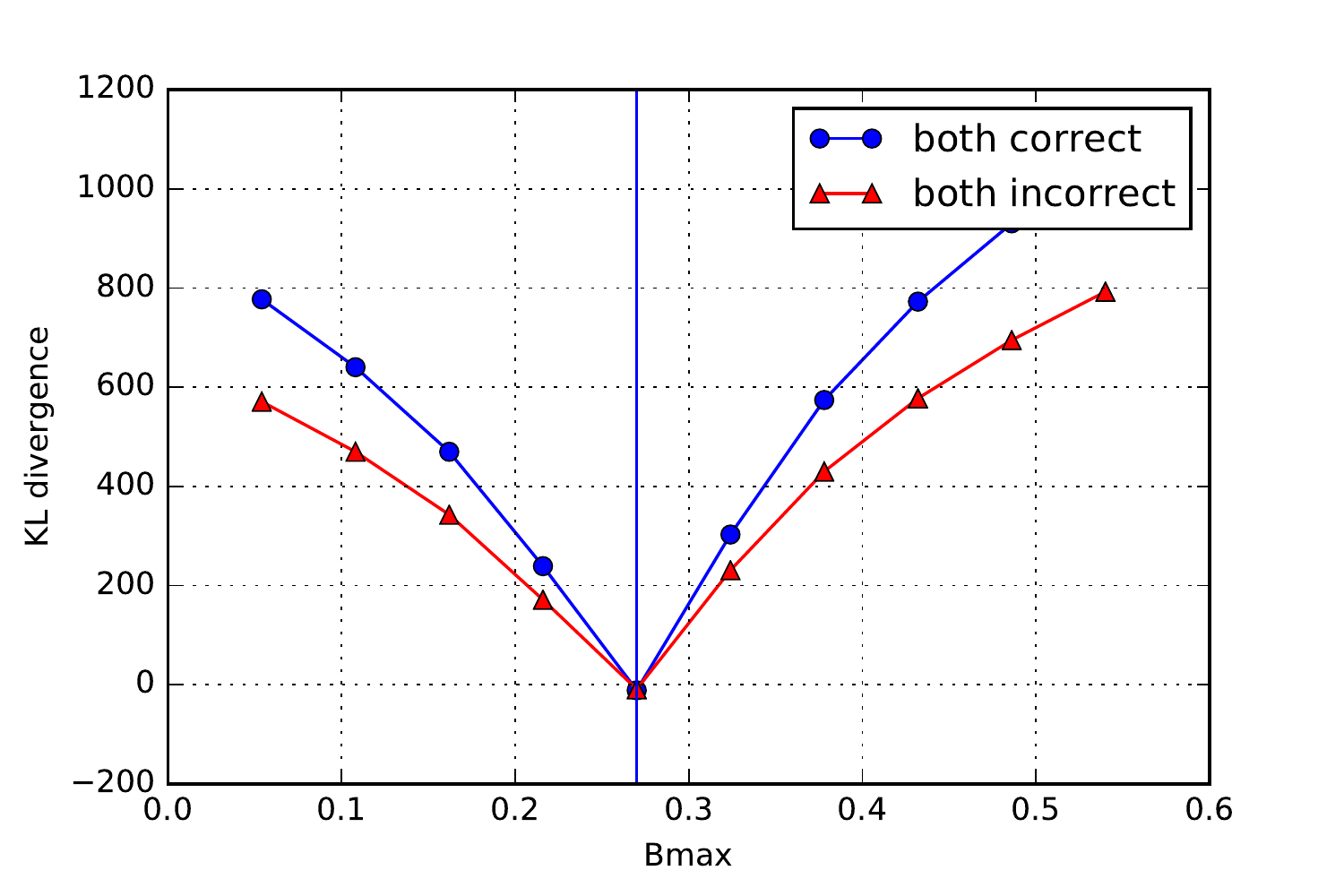} 
   \includegraphics[width=0.48\columnwidth]{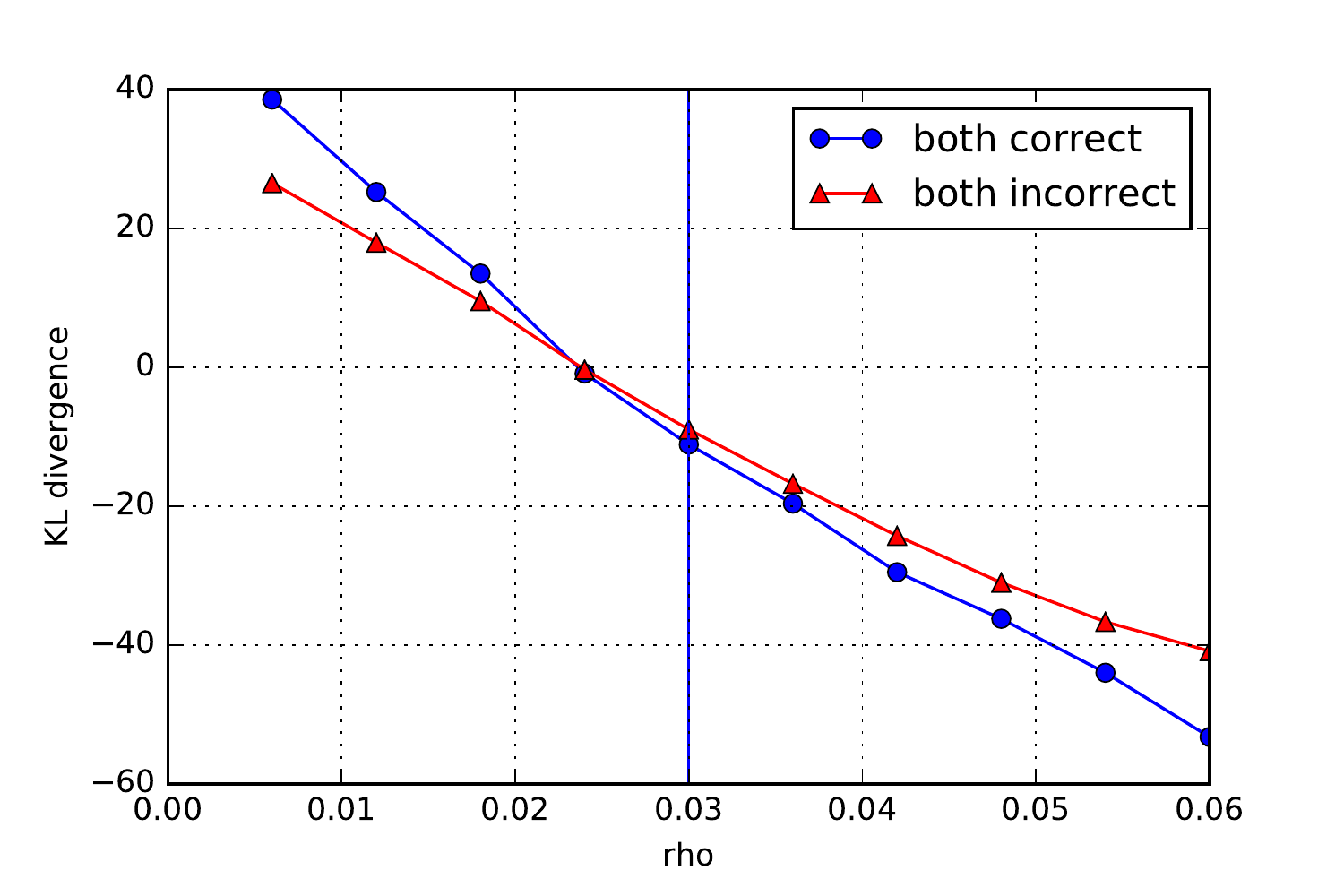} 
   \caption{KL divergence of frequency distributions of resource from the fitted model with that of the original hidden resources.  True parameters are given by the vertical lines.}
   \label{fig:kldiv:min}
\end{figure}

KL divergence comparing distribution of resource values to an \textit{a priori} distribution produced a minimum at the true parameter value for most of the parameters.  However, $\rho$ did not have its minimum near the true value, but at a greater value.  Similarly to the mutual information analysis, $k_b$ had an optimum value for a range of values with the true value near the lower end (Fig. \ref{fig:kldiv:min}).  

\begin{figure}[htbp]
   \centering
   \includegraphics[width=0.48\columnwidth]{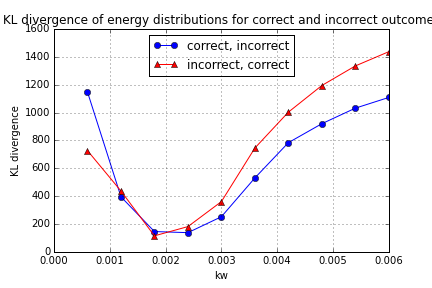} 
   \includegraphics[width=0.48\columnwidth]{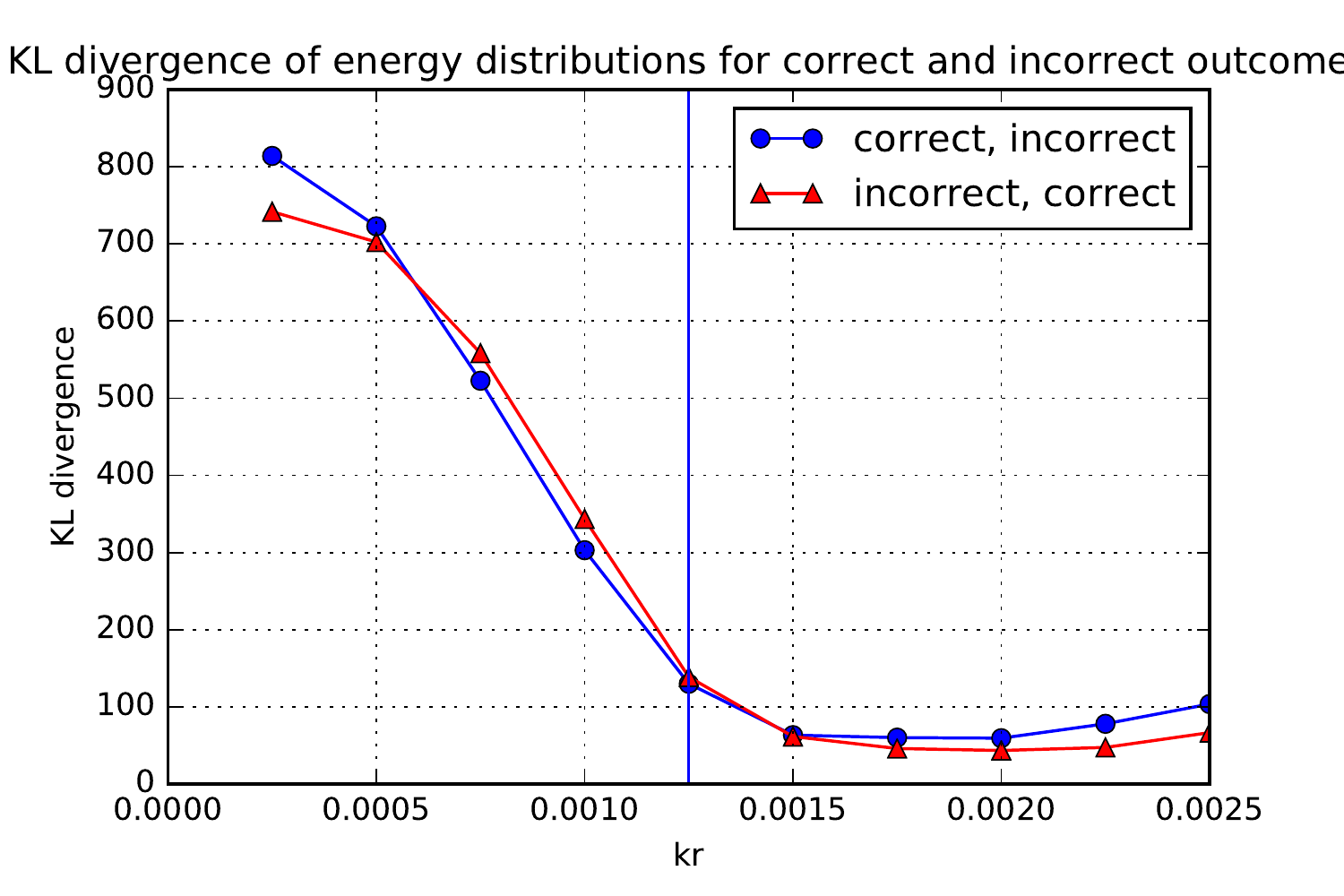} 
   \includegraphics[width=0.48\columnwidth]{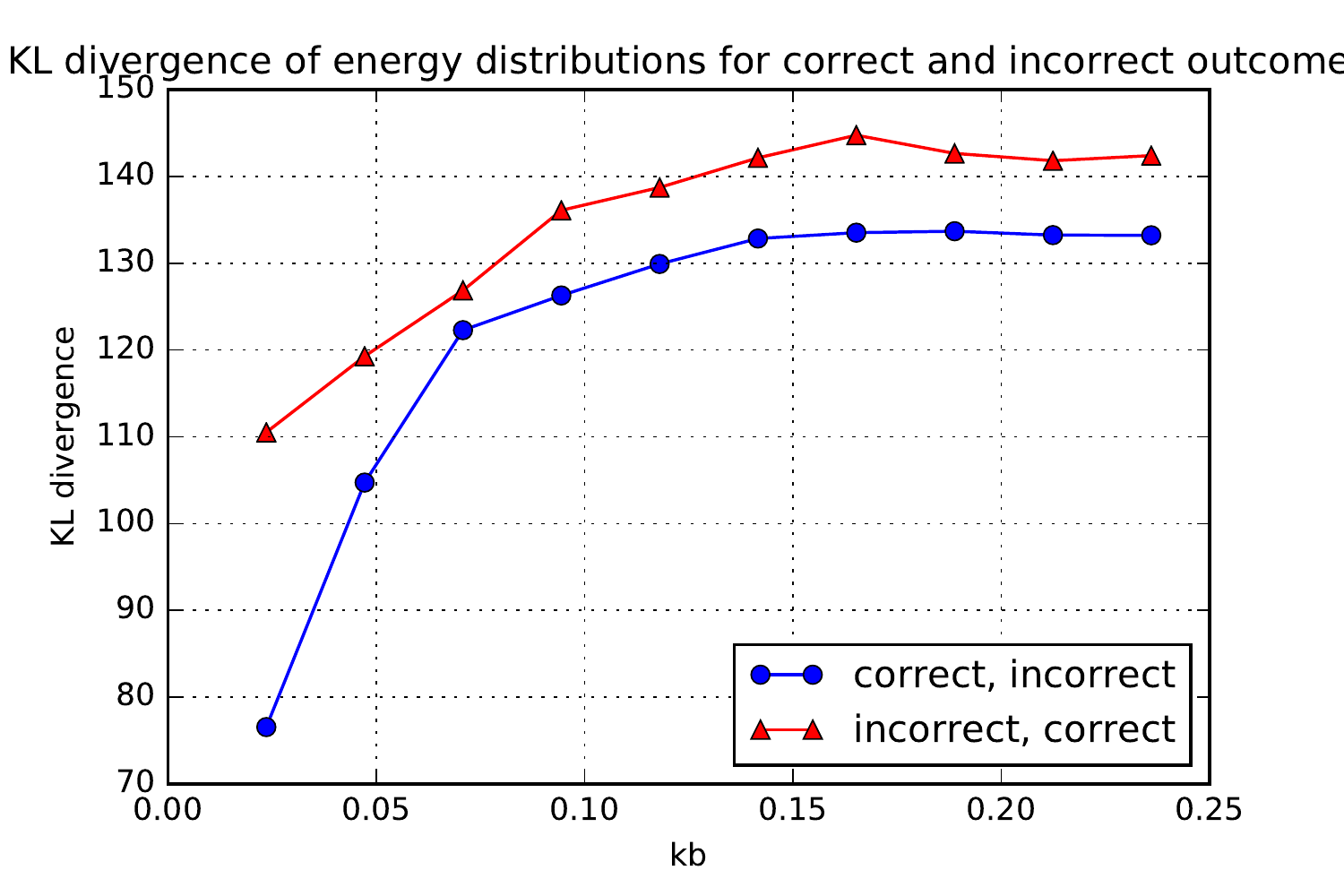} 
   \includegraphics[width=0.48\columnwidth]{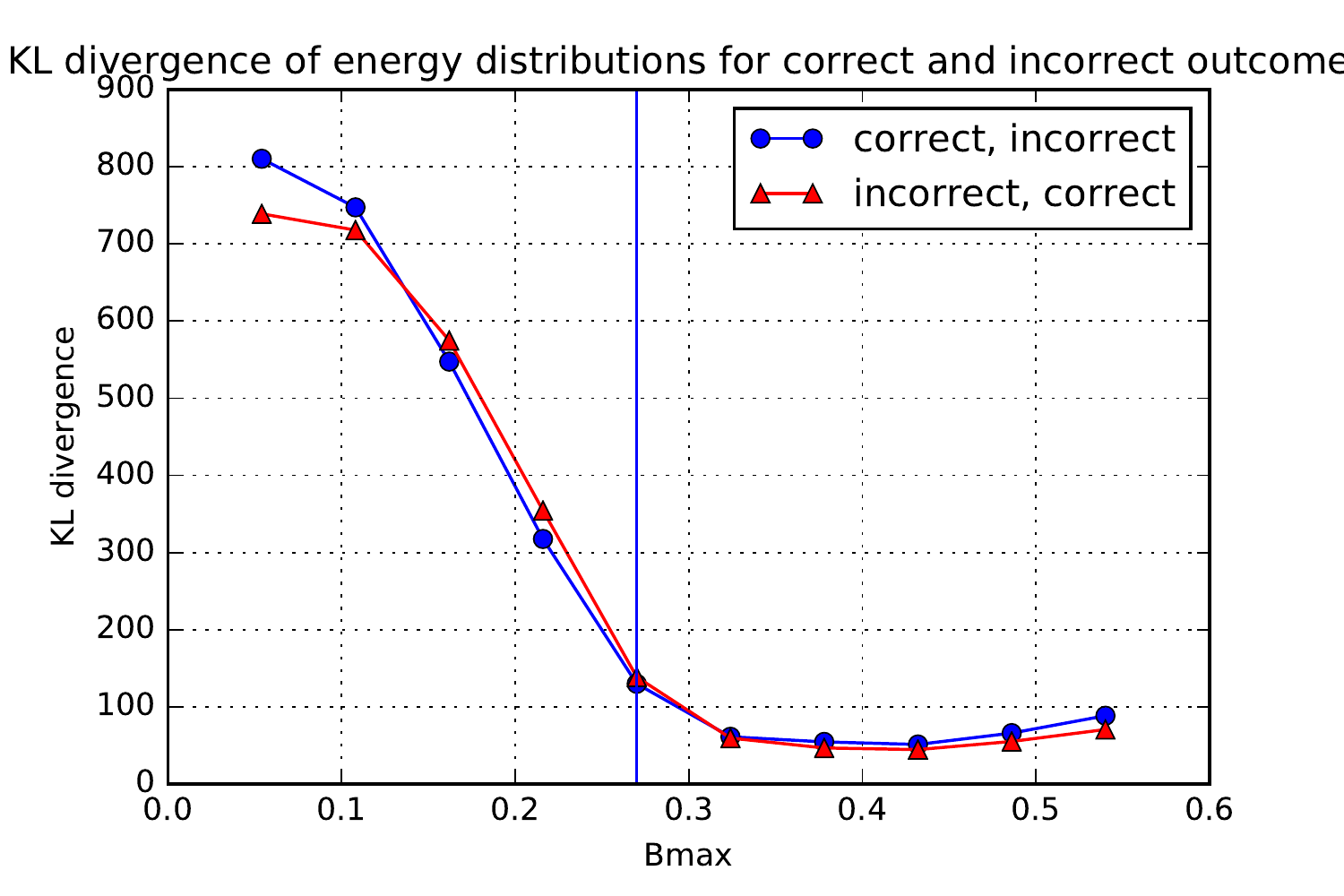} 
   \includegraphics[width=0.48\columnwidth]{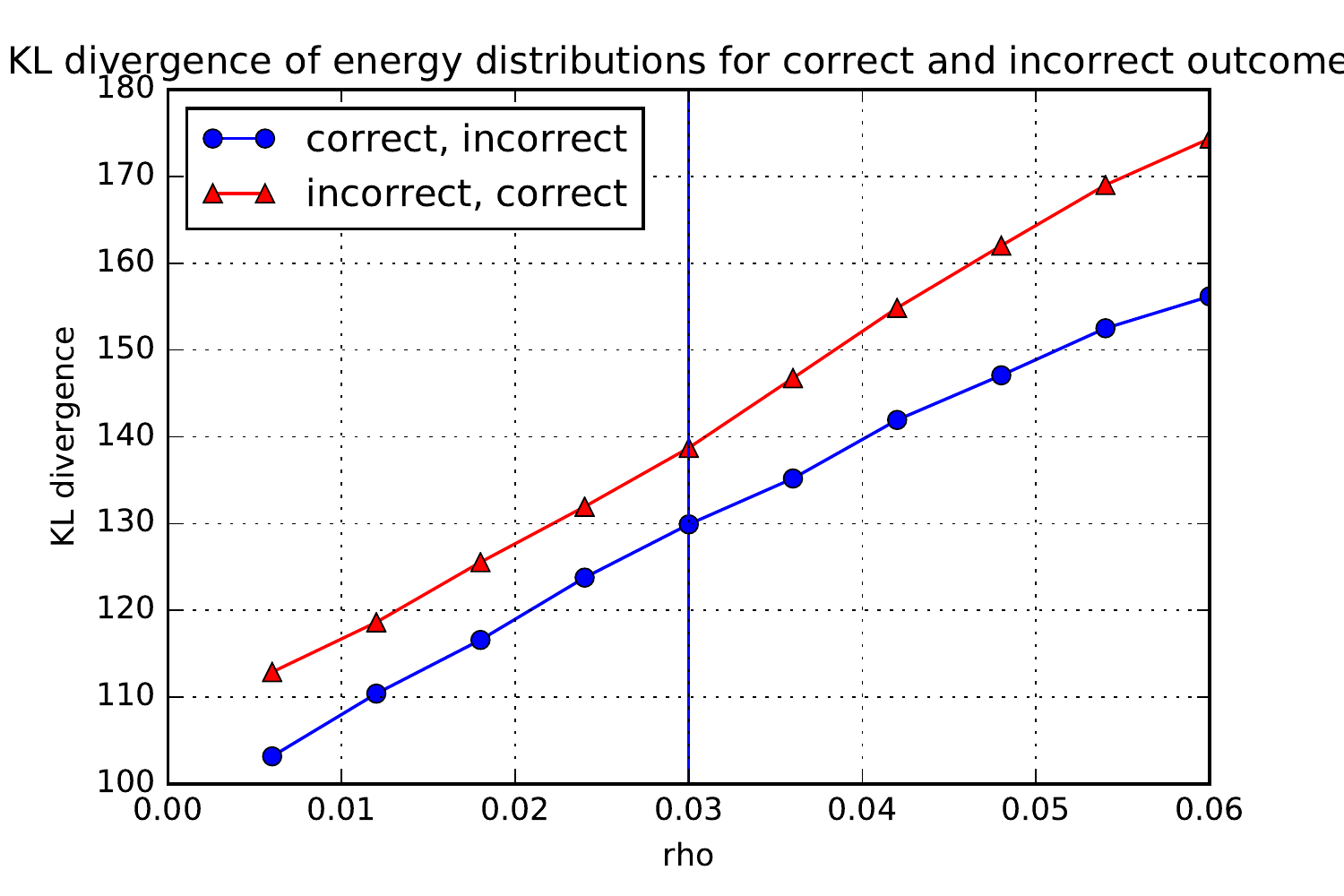} 
   \caption{KL divergence of frequency distributions of resource from the fitted model with that of the original hidden resources.  True parameters are given by the vertical lines.}
   \label{fig:kldiv:max}
\end{figure}

KL divergence comparing the distributions of resource values given outcome fails to produce a maximum at the true parameter; in fact, it produces a local minimum near the true value for $k_w$, $k_r$, and $B_{max}$.  The parameter $k_b$ was the only parameter to produce a maximum at the true value, but similarly to its behavior with other objective functions, it also produced optima at larger values (Fig. \ref{fig:kldiv:max}).

While tests with KL divergence comparing the model to the original resource values produced the clearest optimum at the true parameter values, it poses a problem that the distribution of resource values at the end of each task must be known \textit{a priori}.  In the case of the generated data, the distribution of resource values was found to be Gaussian, allowing the divergence to be approximated using an approximation for the true distribution; however, even if the functional form of the true distribution is known \textit{a priori}, its parameters may not be.  On the other hand, using KL divergence as an objective function to encourage disjointness of resource probability distributions given outcome failed to recover the underlying parameters of the model.    It is worth noting that both successful fitting methods have poor sensitivity for fitting $\rho$.  The unusual behavior in the fitting of $\rho$ and $k_b$ seems to indicate that they are less important to the overall fit of the model, except when $k_b$ is too low.  This may indicate that establishing a lower bound may be more important when fitting this parameter.

\section{Conclusion}
In this work we investigated fitting the hidden layer of deep nonlinear models using mutual information and KL divergence.  Although this was a preliminary investigation, we demonstrated the potential value of using mutual information to fit deep nonlinear models. Although it may be sensitive to biases inherent in the underlying mutual information estimators, the property that mutual information is invariant under isomorphic transformations allows it to estimate the sensitivity of an observation to a model when there is an unknown hidden layer between them. When using KL divergence to compare a model to prior-known statistics of the hidden variables, true parameter values are similarly recovered.  However, using KL divergence to maximize explainability by maximizing the divergence between frequency distributions of the hidden layer given the observed layer completely fails in recovering the original parameter values.

Future work will further explore the formal relationship between mutual information and variations in parameters in the initial, hidden, and observation layers of deep nonlinear models.  The present work provides a promising proof of principle in simple and complex models.  The applications for this work include cognitive science, system identification, and machine learning.

\section{Acknowledgments}
This work was funded by the Analysis in Motion Initiative at Pacific Northwest National Laboratory.

{
\bibliographystyle{plain}
\bibliography{references}
}

\end{document}